%%%                 Reference number:  ??????????
%%%              From R. Cipolatti, I-S. Liu and M. Rincon
%%% 
%%%                    This a Plain LaTeX file
%%%
%%% *************************************************************
%%%                        BEGIN PREAMBLE
%%% *************************************************************

%  ---------------------------------
%  Parameters and fonts
%  ---------------------------------

\baselineskip=14 true pt plus .5pt minus .5pt
\parskip=4pt plus .5pt minus .5pt
\parindent=15pt

\input epsf.tex

\font\bigbf=cmbx10 scaled \magstephalf     % titulo das secoes
\font\Bigbf=cmbx12 scaled \magstep1        % titulo dos capitulos
\font\letranota=cmr8

\font\TenVect=cmbxti10
\font\SevenVect=cmbxti7
\newfam\Vectfam
\def\Vect{\fam\Vectfam\TenVect}
\textfont\Vectfam=\TenVect
\scriptfont\Vectfam=\SevenVect
\scriptscriptfont\Vectfam=\SevenVect

%  ---------------------------------
%  Macros for numeric sets
%  ---------------------------------

\font\TenEns=msbm10
\font\SevenEns=msbm7
\font\FiveEns=msbm5
\newfam\Ensfam
\def\Ens{\fam\Ensfam\TenEns}
\textfont\Ensfam=\TenEns
\scriptfont\Ensfam=\SevenEns
\scriptscriptfont\Ensfam=\FiveEns

\def\R{{\Ens R}}

\def\N{{\Ens N}}

%  ---------------------------------
%  Some little macros  
%  ---------------------------------

\let\Fi=\varphi
\let\eps=\varepsilon
\def\mod#1{\vert #1\vert}
\def\Mod#1{\left|#1\right|}
\def\norma#1#2{\Vert #1\Vert_{#2}}
\setbox111=\hbox to 2truemm{\hrulefill}
\setbox222=\hbox to 2truemm{\vrule height 2truemm width .4truept\hfil\vrule height 2truemm width .4truept}
\def\cqd{\vbox{\offinterlineskip\copy111\copy222\copy111}}
\let\vetor=\letravetor

\def\div{\mathop{\hbox{\rm div}}\nolimits}
\def\tra{\mathop{\rm tr}\nolimits}\let\tr=\tra
\def\vecdiv{\mathop{\hbox{\bf div}}\nolimits}
\def\trepar#1^#2{\mathrel{\mathop{\kern 0pt#1}\limits^{#2}}}
\def\stl{\scriptstyle}
\def\meas{\mathop{\hbox{\rm meas}}\nolimits}
\def\dif{\hbox{\rm d}}
\def\Khi{\raise 2pt\hbox{$\chi$}}

%  ---------------------------------
%  Macros for automatic numbering of formulas and theorems 
%  ---------------------------------

\newcount\numerosection
\newcount\eqnumer
\newcount\lemnumer
\newcount\fignumer
\numerosection=0
\eqnumer=0
\lemnumer=0
\fignumer=0
\def\numsection{\global\advance\numerosection by1
\global\eqnumer=0
\global\lemnumer=0
\global\fignumer=0
\the\numerosection}
\def\strutdepth{\dp\strutbox}
\def\marginalsigne#1{\strut
    \vadjust{\kern-\strutdepth\specialsigne{#1}}}
\def\specialsigne#1{\vtop to \strutdepth{
\baselineskip\strutdepth\vss\llap{#1 }\null}}
\font\margefont=cmr10 at 6pt
\newif\ifshowingMacros
\showingMacrosfalse

\def\cite#1{\csname#1\endcsname}
\def\label#1{\gdef\currentlabel{#1}}
\def\Lefteqlabel#1{\global\advance\eqnumer by 1
\label{#1}
\ifx\currentlabel\relax
\else
\expandafter\xdef
\csname\currentlabel\endcsname{(\the\numerosection.\the\eqnumer)}
\fi
\global\let\currentlabel\relax
\ifshowingMacros
\leqno\llap{%
\margefont #1\hphantom{M}}(\the\numerosection.\the\eqnumer)%
\else
\leqno(\the\numerosection.\the\eqnumer)%
\fi
}

\def\Righteqlabel#1{\global\advance\eqnumer by 1
\label{#1}
\ifx\currentlabel\relax
\else
\expandafter\xdef
\csname\currentlabel\endcsname{(\the\numerosection.\the\eqnumer)}%
\fi
\global\let\currentlabel\relax
\ifshowingMacros
\eqno(\the\numerosection.\the\eqnumer)%
\rlap{\margefont\hphantom{M}#1}
\else
\eqno(\the\numerosection.\the\eqnumer)%
\fi
}
\let\reqlabel=\Righteqlabel
\def\numlabel#1{\global\advance\eqnumer by1     
\label{#1}                                     
\ifx\currentlabel\relax                        
\else
\expandafter\xdef
\csname\currentlabel\endcsname{(\the\numerosection.\the\eqnumer)}
\fi
\global\let\currentlabel\relax
\ifshowingMacros
   (\the\numerosection.\the\eqnumer)\rlap{\margefont\hphantom{M}#1}
\else
(\the\numerosection.\the\eqnumer)%
\fi
}
\def\numero{\global\advance\eqnumer by1       
\number\numerosection.\number\eqnumer}   
\def\lemlabel#1{\global\advance\lemnumer by 1
\ifshowingMacros%
   \marginalsigne{\margefont #1}%
\else
    \relax%\null
\fi
\label{#1}%
\ifx\currentlabel\relax%
\else
\expandafter\xdef%
\csname\currentlabel\endcsname{\the\numerosection.\the\lemnumer}%
\fi
\global\let\currentlabel\relax%
\the\numerosection.\the\lemnumer}
\def\numlem{\global\advance\lemnumer by1   
\the\numerosection.\number\lemnumer}  

%  ---------------------------------
%  Macros for automatic numbering of bibliography
%  ---------------------------------

\let\itemBibli=\item
\def\bibl#1#2\endbibl{\par{\itemBibli{#1} #2\par}}
\def\ref.#1.{{\csname#1\endcsname}}   
\newcount\bib   \bib=0
\def\bibmac#1{\advance\bib by 1       
\expandafter
\xdef\csname #1\endcsname{\the\bib}}
\def\BibMac#1{\advance\bib by1
\bibl{[\the\bib]}%
{\csname#1\endcsname}\endbibl}
\def\MakeBibliography#1{
\noindent{\bf #1}
\bigskip
\bib=0
\let\bibmac=\BibMac
\BiblioFil
\BiblioOrd
}

\newif\ifbouquin

\bouquinfalse
\font\tenssi=cmssi10 at 10 true pt
 at 10 true pt
\def\bibliostyle#1#2#3{{\rm #1}\ifbouquin{\tenssi #2\/}\global\bouquinfalse\else{\it #2\/}\fi{\rm #3}}

%  --------------------------
%  The bibliograpy set 
%  ---------------------------

\def\BiblioFil{
\def\GRS{\bibliostyle{A.E.~Green, R.S.~Rivlin, R.T.~Shield, General theory of small elastic deformations, Pro.~Roy.~Soc.~London, Ser.~A, 211, 1952, pp.~128--154.}{}{}}
\def\LCR{\bibliostyle{I-S.~Liu, R.~Cipolatti, M.A.~Rincon, Successive linear approximation for finite elasticity, Computational and Applied Mathematics, Volume 29, N.~3, 2010, pp.~465--478.}{}{}}
\def\LiuArti{\bibliostyle{I-S.~Liu, Successive linear approximation for boundary value problems of nonlinear elasticity in relative-descriptional formulation, Int.~J.~Eng.~Sci., doi:10.1016/j.ijengsci.2011.02.006, 2011.}{}{}}
\def\Ciarlet{\bibliostyle{P.G.~Ciarlet, {\sl Mathematical Elasticity, Volume 1: Three-Dimensional Elasticity\/}, North-Holland, Amsterdam, 1988.}{}{}}
\def\Ogden{\bibliostyle{R.W.~Ogden, {\sl Non-Linear Elastic Deformations\/}, Ellis Horwood, New York, 1984.}{}{}}
\def\LiuBook{\bibliostyle{I-S.~Liu, {\sl Continuum Mechanics\/}, Springer, Berlin Heidelberg, 2002.}{}{}}
\def\LiuNote{\bibliostyle{I-S.~Liu, A note on Mooney-Rivlin material model, (submitted).}{}{}}
\def\TrusNoll{\bibliostyle{C.~Truesdell, W.~Noll, {\sl The Non-Linear Field Theories of Mechanics\/}, third ed., Springer, Berlin, 2004.}{}{}}
\def\LCRP{\bibliostyle{I-S.~Liu, R.~Cipolatti, M.A.~Rincon, L.A.~Palermo, Numerical simulation of salt migration -- Large deformation in viscoelastic solid bodies, (submitted).}{}{}}
\def\Adams{\bibliostyle{R.A.~Adams, {\sl Sobolev Spaces\/}, Academic Press, New York, 1975.}{}{}}
\def\Haupt{\bibliostyle{P.~Haupt, {\sl Continuum Mechanics and Theory of Materials\/}, Second Edition, Springer, 2002.}{}{}}
}

\def\BiblioOrd{
\bibmac{LCR}      
\bibmac{LiuArti}
\bibmac{LCRP}
\bibmac{Haupt}
\bibmac{TrusNoll} 
\bibmac{LiuNote}   
\bibmac{GRS}
\bibmac{Ciarlet} 
\bibmac{Adams}                  
%\bibmac{LiuBook}     
%\bibmac{Ogden}       

}

\BiblioOrd

%%%  ******************************************************
%%%                    BEGIN DOCUMENT 
%%%  ******************************************************

\centerline{\Bigbf Mathematical analysis of successive linear approximation for}
\centerline{\Bigbf Mooney-Rivlin material model in finite elasticity}

\bigskip
\centerline{Rolci Cipolatti\footnote{*}{\letranota Corresponding author.}, I-Shih Liu, Mauro A.~Rincon}
\bigskip

\footnote{}{\letranota E-mail address: cipolatti@im.ufrj.br (R.~Cipolatti).}

{\baselineskip=14truept
\letranota
\centerline{Instituto de Matem\'atica, Universidade Federal do Rio de Janeiro}
\centerline{C.P. 68530, Rio de Janeiro, 21945-970, RJ, Brazil}
}

\bigskip
\hrule
\bigskip
{%\parindent=1truecm \narrower
\letranota\baselineskip=12truept
\noindent{\bf Abstract}

For calculating large deformations in finite elasticity, we have proposed a method of successive
linear approximation,  by considering the relative descriptional formulation.
In this article we briefly describe this method and we prove the existence and uniqueness of weak solutions for
boundary value problems for nearly incompressible Mooney-Rivlin materials, that arise in each step of the method.

\medskip
\noindent{\letranota Key words: Linearized constitutive equations, Mooney-Rivlin material, Relative descriptional formulation, Existence and uniqueness.}

}
\bigskip
\hrule
\bigskip

\noindent{\bigbf\numsection.\ Introduction}\par
\bigskip

\noindent The constitutive equation of a solid body is usually expressed relative to a preferred reference
configuration which exhibits specific material symmetries such as isotropy. The constitutive
functions are in general nonlinear and linearizations can be used as valid approximation only for small deformations.
Therefore, the problem for large deformations leads to boundary value problems involving systems of nonlinear partial
differential equations.

In order to circumvent the difficulties due to the nonlinearities, we have proposed a new method
for solving numerically the boundary value problem for large deformations. It is based on a successive
linear approximation by considering the relative descriptional formulation.
Roughly speaking, the constitutive equations are calculated at each state which will be regarded as the reference
configuration for the next state, and assuming that the deformation to the next state is small, the updated constitutive
equations can be linearized.

As examples for the proposed method, numerical simulations were done (see [\cite{LCR}], [\cite{LiuArti}]) for two classical
problems concerning Mooney-Rivlin materials, for which the exact solutions are known, namely, the pure shear of a square
and the bending of a retangular block into a circular section. The comparison of the numerical results with the exact
solutions of these two examples confirms the efficiency of our method.

In the present paper we consider the mathematical analysis of the boundary value problem
obtained by linearizing the constitutive equations of nearly incompressible Mooney-Rivlin materials
relative to the present configuration and prove the existence and uniqueness of weak solutions.

We organize this paper as follows. In Section~2 we introduce briefly the notion of relative description and we formally
deduce in Section~3 the linearization of the constitutive function of a nearly incompressible Mooney-Rivlin material.
In Section~4 we consider a boundary value problem involving a system of partial differential equations, that are
obtained by linearizing the constitutive equations of a nearly incompressible Mooney-Rivlin material, and which corresponds to
one of the steps of the successive linear approximation method. The main result of this paper is contained in section 5,
where we prove the existence and uniqueness of weak solutions of this boundary value problem, by considering its variational
formulation. For simplicity, we restrict our analysis to the two-dimensional case, but the arguments  presented
can be extended to three dimensions.

%======================================================================================================
\null\goodbreak
\bigskip
\noindent{\bigbf\numsection.\ Relative description and successive linear approximation}\par
\bigskip
In this section we introduce the notion of {\sl relative description\/}, we formally obtain
the linearization of a general constitutive equation of a solid body respect to this configuration
and we describe the successive linear approximation method.

Let $\kappa_0$ be a reference configuration of a solid body $\cal B$, ${\cal B}_0=\kappa_0({\cal B})$, and let
$$\vetor x=\Khi(X,t),\quad X\in{\cal B}_0$$
be the parametrization of its deformation. Let $\kappa_t$ be the deformed configuration at time $t$
(which we shall always refer as the present time),
${\cal B}_t=\kappa_t({\cal B})$, and
$$F(X,t)=\nabla_X\Khi(X,t)$$
be the deformation gradient with respect to the configuration $\kappa_0$.

Let $\kappa_\tau$ be the deformed configuration at time $\tau>t$. We define the {\sl relative deformation\/} from $\kappa_t$
to $\kappa_\tau$ as the function $\Khi_t:{\cal B}_t\rightarrow{\cal B}_\tau$ given by
$$\Khi_t(\vetor x,\tau):=\Khi(X,\tau),\quad \vetor x\in{\cal B}_t\reqlabel{*1}$$
and the corresponding {\sl relative displacement\/} as
$$\vetor u_t(\vetor x,\tau):=\Khi_t(\vetor x,\tau)-\vetor x.\reqlabel{*2}$$

Taking the gradient relative to $\vetor x$ in both sides of \cite{*2}, we obtain
$$H_t(\vetor x,\tau)= F_t(\vetor x,\tau)-I,\reqlabel{*3}$$
where $I$ is the identity tensor and
$$H_t(\vetor x,\tau):=\nabla_{\vetor x}{\vetor u}_t(\vetor x,\tau),\quad F_t(\vetor x,\tau):=\nabla_{\vetor x}\Khi_t(\vetor x,\tau)$$
are called the {\sl displacement gradient\/} and the {\sl deformation gradient\/} in the relative description, relative to the present configuration.

On the other hand, taking the gradient relative to $X$ in both sides of \cite{*2}, we obtain from \cite{*1} and the chain rule,
$$H_t(\vetor x,\tau)F(X,t)=F(X,\tau)-F(X,t),$$
from which we get easily
$$F(X,\tau)=\bigl(I+H_t(\vetor x,\tau)\bigr)F(X,t).\reqlabel{*5}$$

We can represent this situation by the following picture:
\setbox1=\hbox{\epsfbox{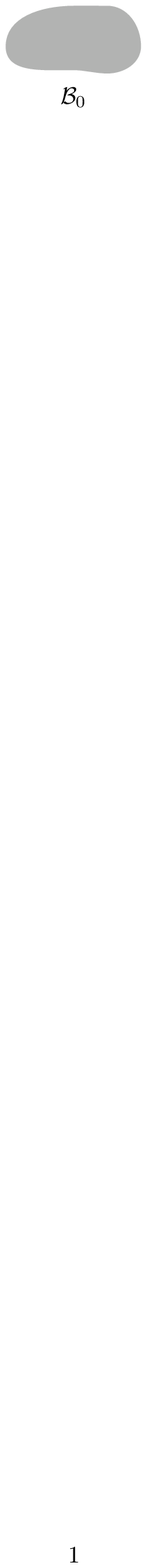}}
\setbox2=\hbox{\epsfbox{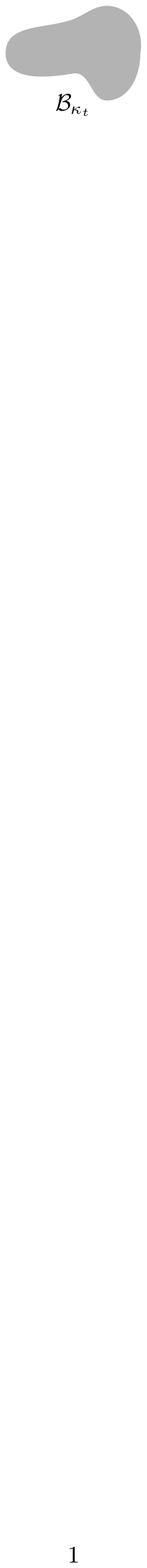}}
\setbox3=\hbox{\epsfbox{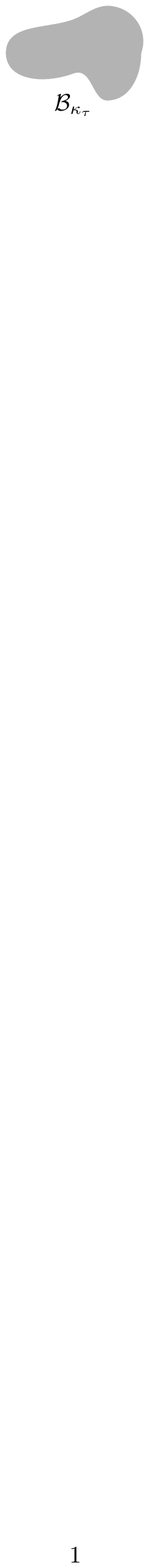}}
$$\hbox{\box1\qquad\raise 7mm
    \hbox{$\trepar{\strut\longrightarrow}_{}^{F(t)}$}
	\qquad\box2\qquad\raise 7mm
    \hbox{$\trepar{\strut\longrightarrow}_{}^{I+H(\tau)}$}\qquad\box3}$$

\noindent where and hereafter, for simplicity, sometimes only the time dependence is indicated. Position dependence is usually
self-evident and will be indicated for clarity only if necessary.

By considering the time $\tau=t+\Delta t$ for small enough $\Delta t$, we can assume that the relative displacement gradient is small,
$$H(\tau):=H_t(\vetor x,\tau),\quad \norma{H(\tau)}{}\ll 1$$

Let $T$ be the Cauchy stress tensor given by the constitutive equation
$$T={\cal F}_{\kappa_0}(F).\reqlabel{2.4}$$
Assuming that the operator ${\cal F}_{\kappa_0}$ is differentiable, we can calculate the linearization of the constitutive equation \cite{2.4}
relative to the current configuration $\kappa_t$, and assuming that $\Delta t$ is small enough, we have formally
$$T(\tau)\approx T(t)+\dif{\cal F}_{\kappa_0}(F(t))[F(\tau)-F(t)]=T(t)+\dif{\cal F}_{\kappa_0}(F(t))[H(\tau)F(t)],\reqlabel{2.5}$$
where $\dif{\cal F}_{\kappa_0}(F)$ denotes the Fr\'echet-differential of ${\cal F}_{\kappa_0}$ calculated at $F$.
For convenience, we shall write \cite{2.5} as
$$T(\tau)=T(t)+L(F(t))[H(\tau)],\reqlabel{2.5meio}$$
where $L(F)[H]:=\hbox{\rm d}{\cal F}_{\kappa_0}(F)[HF]$ defines a fourth order elasticity tensor $L(F)$ relative to the current configuration $\kappa_t$.

The successive linear approximation method is the discrete construction of the parametrization $\Khi(X, t)$ based on the previous arguments.
More precisely, let $t_0<\cdots<t_{n-1}<t_n<t_{n+1}<\cdots$ be a sequence of steps with small enough constant spacing $\Delta t$, where at which step we set
$t=t_n$ and $\tau=t_{n+1}$. Let the deformation gradient $F(X,t_n)$ and the elastic Cauchy stress tensor $T(X,t_n)$,
relative to the preferred configuration $\kappa_{t_0}$, assumed to be known. If in any way we calculate the relative displacement
$\vetor u_{t_n}(\vetor x,t_{n+1})$, $\vetor x\in\kappa_{t_n}(\cal B)$,
it allows us the update the new reference configuration $\kappa_{t_{n+1}}$ relative to the next step by
using \cite{*1} and \cite{*2}, i.e.,
$$\Khi(X,t_{n+1}):=\vetor u_{t_n}(\vetor x,t_{n+1})+\vetor x, \quad\vetor x\in\kappa_{t_n}(\cal B),$$
while the deformation gradient \cite{*5} and the Cauchy stress \cite{2.5meio}, relative to the preferred configuration $\kappa_{t_0}$,
can be determined at instant $t_{n+1}$ respectively by
$$\eqalign{
  F(X,t_{n+1}) & :=\bigl(I+H_{t_n}(\vetor x,t_{n+1})\bigr)F(X,t_n),\cr
  T(X,t_{n+1}) & :=T(X,t_n)+L(F(t_n))[H_{t_n}(\vetor x,t_{n+1})].\cr}$$
Therefore, after updating the boundary data and the eventual body forces acting on the body, we repeat the cycle from the updated
reference configuration $\kappa_{t_{n+1}}$.

We remark that this method can easily be extended to constitutive equation $T={\cal F}(F,\dot F)$ for viscoelastic solid bodies in general [\cite{LCRP}].

%======================================================================================================
\null\goodbreak
\bigskip
\noindent{\bigbf\numsection.\ Application to nearly incompressible Mooney-Rivlin materials}\par
\bigskip

From now on we consider a  Mooney-Rivlin material whose constitutive equation relative
to the preferred reference configuration $\kappa_0$ is given by
$$T={\cal F}_{\kappa_0}(F)=-pI+\widetilde{\cal F}(F),\quad\widetilde{\cal F}(F)=s_1B+s_2B^{-1},$$
where $B=FF^T$ is the left Cauchy-Green strain tensor and the material parameters $s_1$ and $s_2$ are constant
satisfying
$$s_1>0\quad\hbox{\rm and}\quad s_2<s_1.\reqlabel{CondS1S2}$$

\noindent{\bf Remark\ \lemlabel{Obs0}:} It is usually assumed $s_2\le 0< s_1$, the so-called E-inequalities (see [\cite{Haupt}] and [\cite{TrusNoll}]), based on the assumption
that the free energy function is positive definite for any deformation. Liu [\cite{LiuNote}] 
has pointed out that ``any'' deformation is unrealistic from physical point of view, and a thermodynamical stability anaysis only requires $s_2<s_1$. 
Therefore, we shall include the case $0<s_2<s_1$ in our analysis.

\smallskip\goodbreak
A direct calculation of the Fr\'echet-differential of $\widetilde{\cal F}$ at $F$ gives
$$\dif\widetilde{\cal F}(F)[H]=s_1\bigl(HB+BH^{T}\bigr)-s_2\bigl(B^{-1}H+H^TB^{-1}\bigr).$$

For compressible body in general, the pressure $p$ may depend on the deformation gradient $F$. However,
for compressible elastic bodies, we shall assume that the pressure depends only on the determinant of the
deformation gradient or, by the mass balance, depends only on the mass density $\rho$,
$$p=p(\rho),\quad \rho(t)={\rho_0\over \det F(t)},$$
where $\rho_0$ denotes the mass density in the referential configuration $\kappa_0$.

For time $\tau=t+\Delta t$ and from \cite{*5}, we have
$$\eqalign{
  \rho(\tau)-\rho(t) & = \rho_0\bigl(\det F(\tau)^{-1}-\det F(t)^{-1}\bigr)=\rho(t)\bigl(\det F(t)F(\tau)^{-1}-1\bigr)\cr
                     & = \rho(t)\bigl(\det(I+H(\tau))^{-1}-1\bigr)=-\rho(t)\tra H(\tau)+o(2), \cr
}$$
where  $\tr H$ means the trace of $H$ and $o(2)$ denotes higher order terms in the small displacement gradient $H(\tau)$.
Therefore, assuming that $p$ is differentiable as function of $\rho$, we have
$$p(\tau)-p(t)=\left({dp\over d\rho}\right)_t\bigl(\rho(\tau)-\rho(t)\bigr)+o(2)=-\left(\rho{dp\over d\rho}\right)_t\tra H(\tau)+o(2),$$
or
$$p(\tau)=p(t)-\beta(t)\tra H(\tau)+o(2),$$
where $\beta(t)=\rho(t)(dp/d\rho)_t$ is a material parameter depending on the mass density $\rho$.

A body is called {\sl nearly incompressible\/} if its density is nearly insensitive to change of pressure. Hence, if
we regard the density as a function of pressure, $\rho=\rho(p)$, then its derivative with respect to the pressure is nearly zero.
This means that, for nearly incompressible materials, the parameter $\beta$ must be large,
$$\beta=\beta(\vetor x,t)\gg 1,\quad\forall\vetor x\in{\cal B}_t.$$

Therefore, the Cauchy stress tensor relative to the current configuration $\kappa_t$ is given by
$$T(\tau)=T(t)+L(F(t))[H(\tau)]+o(2),$$
where
$$L(F)[H]=\beta(\tra H)I+s_1\bigl(HB+BH^{T}\bigr)-s_2\bigl(B^{-1}H+H^TB^{-1}\bigr)$$
and the first Piola-Kirchhoff stress tensor at time $\tau$ relative to the current configuration $\kappa_t$ is given by
$$\eqalign{
T_{\kappa_t}(\tau) & = \det F_t(\tau)T(\tau)F_t(\tau)^{-T}=\det(I+H)T(\tau)(I+H)^{-T}\cr
                   & = \bigl[1+\tra H+o(2)\bigr]\bigl[T(t)+L(F(t))[H]+o(2)\bigr]\bigl[I-H^T+o(2)\bigr]\cr
				   & = T(t)+(\tra H)T(t)-T(t)H^T+L(F(t))[H]+o(2).\cr
}\reqlabel{PiolKirch}$$

%======================================================================================================
\null\goodbreak
\bigskip
\noindent{\bigbf\numsection.\ Linearized boundary value problem and its variational formulation}\par
\bigskip
For simplicity, we denote by $\kappa$ the current configuration $\kappa_t$, $\Omega=\kappa({\cal B})$ be the bounded domain of $\R^3$
representing the interior of the region occupied by the body at current configuration $\kappa$ at the present time $t$, $T_0=T(t)$
and $B_0=B(t)$. Let $\partial\Omega=\Gamma_1\cup\Gamma_2\cup\Gamma_3$, $\vetor n_\kappa$ be the exterior unit normal to $\partial\Omega$
and $\vetor g$ be the gravitational force (per unit mass).

We consider the following boundary value problem for the relative displacement $\vetor u=\vetor u(\vetor x,\tau)$,

$$\left\{\eqalign{
-\vecdiv T_{\kappa}(\tau) & = \rho(\tau)\vetor g\quad \hbox{\rm in\ } \Omega\times\R,\cr
T_{\kappa}(\tau)\vetor{n}_{\kappa} & = \vetor f(\tau)\quad \hbox{\rm on\ } \Gamma_1, \cr
\vetor{u}(\tau)\cdot\vetor{n}_{\kappa} & = 0\quad \hbox{\rm on\ } \Gamma_2,\cr
T_{\kappa}(\tau)\vetor{n}_{\kappa}\times\vetor{n}_{\kappa}& =\vetor 0\quad \hbox{\rm on\ } \Gamma_2,\cr
\vetor{u}(\tau) & = \vetor 0\quad \hbox{\rm on\ } \Gamma_3,\cr
}\right.\reqlabel{ProbCont}$$
where $\vecdiv$ is the divergence operator with respect to $\vetor x$, $T_\kappa(\tau)=T_\kappa(\vetor x,\tau)$ is the Piolla-Kirchhoff stress tensor at time $\tau$ relative to
configuration $\kappa$ at the present time $t$, which, up to linear terms in relative displacement gradient $H=H(\tau)=\nabla_{\vetor x}\vetor u(\tau)$, is given by (see \cite{PiolKirch})
$$T_\kappa = T_0+(\tra H)(T_0+\beta I)-T_0H^T+s_1\bigl(HB_0+B_0H^{T}\bigr)-s_2\bigl(B_0^{-1}H+H^TB_0^{-1}\bigr),$$
and $\vetor f(\tau)$ is the surface traction (per unit surface area).

\smallskip
\noindent{\bf Remark\ \lemlabel{Obs1}:} At every time step, the idea of formulating the boundary value problem in the
form \cite{ProbCont} is similar to the theory of small deformations superposed on finite
deformations (see [\cite{GRS}], [\cite{TrusNoll}]). In this manner, either we are interested in the
evolution of solutions with gradually changing boundary conditions resulting
in large deformation, or, we can treat the boundary values of finite
elasticity as the final value of a successive small incremental boundary
values at each time step (see [\cite{Ciarlet}]).

\smallskip
The boundary value problem \cite{ProbCont} can be formulated as a variational problem. Indeed, let $\Omega$ be a smooth enough bounded domain
in $\R^3$ and define the space
$${\cal V}=\bigl\{\vetor u\in(H^1(\Omega))^3\,;\, \vetor u\cdot\vetor n_\kappa=0 \hbox{\ on\ }
       \Gamma_2\,\,\,\hbox{\rm and}\,\,\, \vetor u=\vetor 0 \hbox{\ on\ } \Gamma_3 \bigr\}.$$
Taking formally the inner product of both sides of the equation in \cite{ProbCont} by $\vetor w\in{\cal V}$ and integrating over $\Omega$,
we obtain after integration by parts,
$$\int_\Omega \tra\bigl(K[H]W^T\bigr)\,d\vetor x=\int_{\Gamma_1}\vetor f(\tau)\cdot\vetor w\,d\Gamma-\int_\Omega\tra\bigl(T_0W^T\bigr)\,d\vetor x,$$
where we are denoting $H=\nabla_{\vetor x}\vetor u$, $W=\nabla_{\vetor x}\vetor w$ and $K[H]$ is given by
$$K[H]:=(\tra H)(T_0+\beta I)-T_0H^T+s_1\bigl(HB_0+B_0H^{T}\bigr)-s_2\bigl(B_0^{-1}H+H^TB_0^{-1}\bigr).$$	
Therefore, for  $\vetor u,\vetor w\in{\cal V}$ we consider respectively the bilinear and the linear forms:
$$\eqalign{
  {\cal L}(\vetor u,\vetor w) & :=\int_\Omega\tra(K[H]W^T)\,d\vetor x,\cr
  {\cal N}(\vetor w) & :=\int_{\Gamma_1}\vetor f(\tau)\cdot\vetor w\,d\Gamma-\int_\Omega\tra\bigl(T_0W^T\bigr)\,d\vetor x+\int_\Omega\rho(\tau)\vetor g\cdot\vetor w\,d\vetor x.\cr
}\reqlabel{Formas}$$
We notice that the forms ${\cal L}$ and ${\cal N}$ can be written in terms of coordinates by
$$\eqalign{
{\cal L}(\vetor u,\vetor w) & = \int_\Omega{\partial u_k\over \partial x_k}\Bigl([{T}_0]_{ij}+\beta\delta_{ij}\Bigr){\partial w_i\over \partial x_j}\,dV
	    -\int_\Omega[{T}_0]_{ik}{\partial u_j\over\partial x_k}{\partial w_i\over\partial x_j}\,dV \cr
	&\quad{}+s_1\int_\Omega\left({\partial u_i\over\partial x_k}[B_0]_{kj}+
	          [B_0]_{ik}{\partial u_j\over\partial x_k}\right){\partial w_i\over\partial x_j}\,dV\cr
	&\qquad{}- s_2\int_\Omega\left([B_0^{-1}]_{ik}{\partial u_k\over\partial x_j}+
	          {\partial u_k\over\partial x_i}[B_0^{-1}]_{kj}\right){\partial w_i\over\partial x_j}\,dV,\cr	
{\cal N}(\vetor w) & =\int_{\Gamma_1}f_iw_i\, d\Gamma-\int_\Omega [ {T}_0]_{ij}{\partial w_i\over\partial x_j}\, d\vetor x+\int_\Omega\rho g_iw_i\,d\vetor x,\cr
 }$$
where in the above formulas we have used  the standard summation convention for repeated indices.

Then, the variational problem is to find the solution $\vetor u\in{\cal V}$ such that
$${\cal L}(\vetor u,\vetor w)={\cal N}(\vetor w),\quad \forall \vetor w\in{\cal V}.\reqlabel{Variacional}$$

In order to prove that the solutions of \cite{Variacional} is a weak solution of \cite{ProbCont}, the following result concerning existence of
a {\sl normal trace\/} is useful.

\smallskip
\noindent{\bf Lemma\ \lemlabel{Lemma0}:} {\sl Let $\Omega\subset\R^N$ be a bounded Lipschitz domain. If
$\vetor F\in L^2(\Omega)^N$ satisfies $\div\vetor F\in L^2(\Omega)$, then $\vetor F\cdot\vetor n_\kappa$ can be defined as an element of $H^{-1/2}(\partial\Omega)$
and there exists a constant $C_1>0$ depending only on $\Omega$ such that
$$\norma{\vetor F\cdot\vetor n_\kappa}{H^{-1/2}}\le C_1\Bigl(\norma{\vetor F}{2}+\norma{\div\vetor F}{2}\Bigr).$$
}

\noindent{\bf Proof:} Assume that $\vetor F\in C^1(\Omega)\cap C^0(\overline\Omega)$. Then, using integration by parts, por any $\psi\in C^1(\Omega)\cap C^0(\overline\Omega)$ we have
$$\int_\Omega \vetor F(\vetor x)\cdot\nabla\psi(\vetor x)\,d\vetor x+\int_\Omega\div\vetor F(\vetor x)\psi(\vetor x)\,d\vetor x=
    \int_{\partial\Omega}\psi(\vetor x)\vetor F(\vetor x)\cdot\vetor n_\kappa(\vetor x)\,d\Gamma.$$ 	
Therefore, denoting the right hand side of the above identity as $\langle\vetor F\cdot\vetor n_\kappa\,;\, \gamma_0(\psi)\rangle$, with the brackets meaning the duality between
$H^{-1/2}(\partial\Omega)$ and $H^{1/2}(\partial\Omega)$ and $\gamma_0:H^1(\Omega)\rightarrow H^{1/2}(\partial\Omega)$ being the trace operator, we have
$$\mod{\langle\vetor F\cdot\vetor n_\kappa\,;\, \gamma_0(\psi)\rangle}\le \norma{\vetor F}{2}\norma{\nabla\psi}{2}+\norma{\div\vetor F}{2}\norma{\psi}{2}.$$
It is well-known that, for a given $\Fi\in H^{1/2}(\partial\Omega)$ we may choose $\psi\in H^1(\Omega)$ such that $\gamma_0(\psi)=\Fi$ and such that
$\norma{\psi}{H^1}\le C_1\norma{\Fi}{H^{1/2}}$, where the constant $C_1$ depends only on $\Omega$. Hence,
$$\mod{\langle\vetor F\cdot\vetor n_\kappa\,;\, \Fi\rangle}\le C_1\bigl(\norma{\vetor F}{2}+\norma{\div\vetor F}{2}\bigr)\norma{\Fi}{H^{1/2}}.\reqlabel{IneqKav}$$	
This means that
$$\norma{\vetor F\cdot\vetor n_\kappa}{H^{-1/2}}\le C_1\bigl(\norma{\vetor F}{2}+\norma{\div\vetor F}{2}\bigr).$$
When $\vetor F$ is no longer in $C^1(\Omega)\cap C^0(\overline\Omega)$, using a density argument (see [\cite{Adams}]), 
we can find a sequence $\{\vetor F_n\}_{n\in\N}$ in $C^1(\Omega)\cap C^0(\overline\Omega)$ such that
$$\vetor F_n\rightarrow \vetor F\quad\hbox{\rm in}\quad L^2(\Omega)^N,\quad \div\vetor F_n\rightarrow \div\vetor F\quad\hbox{\rm in}\quad L^2(\Omega).$$
Inequality \cite{IneqKav} shows that $\{\vetor F_n\cdot\vetor n_\kappa\}_{n\in\N}$ is a Cauchy sequence in $H^{-1/2}(\partial\Omega)$, whose limit,
which is independent of the particular choice of the sequence $\{\vetor F_n\}_{n\in\N}$, will be denoted by  $\vetor F\cdot\vetor n_\kappa$. This finishes the proof.\quad\cqd

\smallskip
\noindent{\bf Lemma\ \lemlabel{Lemma1}:} {\sl Let $\Omega\subset\R^3$ be a domain of class $C^2$. We assume that
$\beta, p_0\in L^\infty(\Omega)$, $\rho\in L^2(\Omega)$ and $T_0,B_0\in L^\infty(\Omega,M_3(\R))$, where $M_3(\R)$ denotes the set of $3\times3$
real matrices. If $\vetor u$ is a solution of \cite{Variacional}, then $\vetor u$ is a weak solution of \cite{ProbCont}. }

\smallskip
\noindent{\bf Proof:} Let $u\in\cal V$ be a solution of \cite{Variacional}. Then, $H:=\nabla_{\vetor x}\vetor u\in L^2(\Omega,M_3(\R))$, which implies that
$T_\kappa=T_0+K[H]\in L^2(\Omega,M_3(\R))$. Since $C_0^\infty(\Omega)^3\subset {\cal V}$, we have
$$-\int_\Omega\vecdiv T_\kappa\cdot \vetor w\,d\vetor x=\int_\Omega\tr (T_\kappa W^T)\,d\vetor x
     =\int_\Omega \rho\vetor g\cdot\vetor w\,d\vetor x,\quad\forall\vetor w\in C_0^\infty(\Omega)^3,$$
where $W=\nabla_{\vetor x}\vetor w$ and the partial derivatives in $\vecdiv$ are taken in the sense of distributions in $\Omega$. Hence, $\vetor u$ satisfies
$$-\vecdiv T_\kappa=\rho\vetor g$$
in the sense of distributions.	Moreover, since we are assuming that $\rho\in L^2(\Omega)$, it follows from the density of $C_0^\infty(\Omega)$ in $L^2(\Omega)$
that $\vecdiv T_\kappa\in L^2(\Omega)^3$. From Lemma~\cite{Lemma0}, \cite{Variacional} reduces to
$$\int_{\partial\Omega}T_\kappa\vetor n_\kappa\cdot\vetor w\,d\Gamma=\int_{\Gamma_1}\vetor f\cdot\vetor w\, d\Gamma,\reqlabel{N**}$$
where the above surface integral on ${\partial\Omega}$ are taken in the sense of the duality between $H^{-1/2}(\partial\Omega)^3$
and $H^{1/2}(\partial\Omega)^3$.
In particular, for any $\vetor w\in\cal V$ such that $\vetor w=\vetor 0$ on $\Gamma_2$, we have
$$\int_{\Gamma_1}\bigl( T_\kappa\vetor n_\kappa-\vetor f\bigr)\cdot\vetor w\,d\Gamma=0,$$
which gives the $\Gamma_1$-boundary condition in \cite{ProbCont}. So, \cite{N**} reduces to
$$\int_{\Gamma_2}T_\kappa\vetor n_\kappa\cdot\vetor w\,d\Gamma=0,\quad\forall\vetor w\in\cal V.\reqlabel{N***}$$

In order to show that \cite{N***} gives the $\Gamma_2$-boundary condition in \cite{ProbCont}, let $\Fi\in H_0^1(\Omega)$, $\Fi<0$, be the
first eigenfunction of $-\Delta$ and define
$$\vetor w_0(\vetor x):=\nabla_{\vetor x}\Fi(\vetor x)\mod{\nabla_{\vetor x}\Fi(\vetor x)}^{-1},\quad\vetor x\in\Omega.$$
Since $\Omega$ is of class $C^2$, we can extend $\vetor w_0$ to the boundary $\partial\Omega$ and we have
from the maximum principle that $\vetor w_0(\vetor x)=\vetor n_\kappa(\vetor x)$, for almost all $\vetor x\in\partial\Omega$.
Let $\widetilde{\vetor w}\in H^1(\Omega)^3$ be an arbitrary function which vanishes on $\Gamma_3$ and consider $\vetor w=\vetor w_0\times \widetilde{\vetor w}$.
Then, it is clear that $\vetor w\in{\cal V}$, since $\vetor w=\vetor 0$ on $\Gamma_3$ and
$$\vetor w\cdot\vetor n_\kappa\,|_{\Gamma_2}=(\vetor n_\kappa\times\widetilde{\vetor w})\cdot\vetor n_\kappa\,|_{\Gamma_2}
    =-(\vetor n_\kappa\times\vetor n_\kappa)\cdot\widetilde{\vetor w}|_{\Gamma_2}=0.$$
Therefore,  from \cite{N***},
$$0=\int_{\Gamma_2}T_\kappa\vetor n_\kappa\cdot\vetor w\,d\Gamma =
     \int_{\Gamma_2}T_\kappa\vetor n_\kappa\cdot(\vetor n_\kappa\times\widetilde{\vetor w})\,d\Gamma =
	  \int_{\Gamma_2}(T_\kappa\vetor n_\kappa\times\vetor n_k)\cdot \widetilde{\vetor w}\,d\Gamma$$
and the proof is complete.\quad\cqd

\goodbreak
\bigskip
\noindent{\bigbf\numsection.\ Existence and uniqueness of solution in two-dimensions}\par
\bigskip
Let $\Omega$ be a bounded Lipschitz domain of $\R^2$ whose boundary $\partial\Omega=\Gamma_1\cup\Gamma_2\cup\Gamma_3$, with
$\meas(\Gamma_i)\not=0$ for $i=1,2,3$, and consider the space
$${\cal V}=\bigl\{\vetor u=(u_1,u_2)\in H^1(\Omega)^2\,;\,\vetor u\cdot\vetor n_\kappa=0 \hbox{\ on\ }
       \Gamma_2\,\,\,\hbox{\rm and}\,\,\, \vetor u=\vetor 0 \hbox{\ on\ } \Gamma_3 \bigr\}.\reqlabel{DefV2D}$$
For $\vetor u,\vetor w \in {\cal V}$, we introduce
$$\eqalign{
\bigl\langle \vetor u|\vetor v\bigr\rangle & := \int_\Omega\bigl(\nabla u_1(\vetor x)\cdot\nabla v_1(\vetor x)+\nabla u_2(\vetor x)\cdot\nabla v_2(\vetor x)\bigr)\, d\vetor x,\cr
\norma{\vetor u}{\cal V}^2& := \norma{\nabla u_1}{L^2}^2+\norma{\nabla u_2}{L^2}^2,\cr
}\reqlabel{DefNormaV}$$	
where $\norma{\,\,\,}{L^2}$ is the usual $L^2$-norm.	It is well-known that the Poincar\'e inequality holds if ${\meas(\Gamma_3)\not=0}$, i.e., there exists a constante $C$ such that
$$\norma{\vetor u}{\cal V}^2\ge C\norma{\vetor u}{L^2}^2,\quad\forall\vetor u\in{\cal V}.$$
In this case, $\langle\cdot;\cdot\rangle$ and $\norma{\,\,}{\cal V}$ define an inner product and a norm in ${\cal V}$, respectively.

From now on we assume that
$$\rho\in L^2(\Omega),\quad \beta, p_0\in L^\infty(\Omega),\quad B_0\in L^\infty\bigl(\Omega, S_2^+(\R)\bigr),\reqlabel{Hypot}$$
where by $S_2^+(\R)$ we denote the set of all symmetric and positive definite $2\times 2$ matrix, and we set
$$T_0:=-p_0I+s_1B_0+s_2B_0^{-1}.$$
It is clear that the forms ${\cal L}$ and $\cal N$ defined in \cite{Formas} are continuous in $\cal V$.

Recalling that $H$ and $W$ are $2\times 2$ matrix whose entries are given by
$$[H]_{ij}={\partial u_i\over\partial x_j},\quad [W]_{ij}={\partial w_i\over\partial x_j},\quad \vetor u,\vetor w\in{\cal V},$$
the bilinear form ${\cal L}(\vetor u,\vetor w)$ defined in \cite{Formas} can be written as
$${\cal L}(\vetor u,\vetor w)=\int_\Omega {\cal A}({\vetor x};H({\vetor x}),W({\vetor x}))\,d\vetor x,$$
where
$$\eqalign{
{\cal A}({\vetor x};H,W) & := \tra(H)\tra\Bigl[(T_0+\beta I)W^T\Bigr]-\tra(T_0 H^TW^T)\cr
              &\quad{}+s_1\tra\Bigl[(HB_0+B_0H^T)W^T\Bigr]- s_2\tra\Bigl[(B_0^{-1}H+H^TB_0^{-1})W^T\Bigr].\cr
}$$
In particular, for $W=H$ we have
$$\eqalign{
{\cal A}({\vetor x};H,H) & = \tra(H)\tra\Bigl[(T_0+\beta I)H^T\Bigr]-\tra(T_0 H^TH^T)\cr
              &\quad{}+s_1\tra\Bigl[(HB_0+B_0H^T)H^T\Bigr]-s_2\tra\Bigl[(B_0^{-1}H+H^TB_0^{-1})H^T\Bigr].\cr
}\reqlabel{DefFormQuadA}$$
Hence, to prove that ${\cal L}$ is coercive, it is sufficient to show that there exists $\alpha>0$ such that
$${\cal A}({\vetor x};H,H)\ge \alpha\norma{H}{}^2,\,\,\forall\,{\vetor x}\in\Omega,$$
i.e., it suffices to show that the bilinear form ${\cal A}({\vetor x};H,W)$ is uniformly coercive as function of $2\times 2$ matrices.
Furthermore, a direct calculation (see the Appendix) gives that the coercivity of ${\cal A}(\vetor x,H,W)$ is equivalent to the semipositivity of the
matrix $A(\vetor x)-\alpha I$, for all $\vetor x\in\Omega$, with $A(\vetor x)$ given by
$$A(\vetor x)\kern-2pt=\kern-2pt\pmatrix{
\stl \beta+2s_1\gamma_1-2{s_2}\gamma_1^{-1} & \stl \beta+{1\over 2}\tra T_0 & \stl 0 &\stl 0 \cr
\stl \beta+{1\over 2}\tra T_0 & \stl \beta+2s_1\gamma_2-2{s_2}\gamma_2^{-1}    & \stl 0 & \stl 0 \cr
\stl 0 & \stl 0 & \stl  2s_1\tra B_0 - 2{s_2}\tra B_0^{-1}-\tra T_0    & \stl  s_1(\gamma_2-\gamma_1)-{s_2}(\gamma_1^{-1}-\gamma_2^{-1})\cr
\stl 0 & \stl 0 & \stl s_1(\gamma_2-\gamma_1)-{s_2}(\gamma_1^{-1}-\gamma_2^{-1}) & \stl\tra T_0                              \cr
}\reqlabel{MatrizA}$$
where $\gamma_1$ e $\gamma_2$ are the eigenvalues of $B_0$ and $I$ is the $4\times 4$ identity matrix. Therefore, we can also write
$${\cal L}(\vetor u,\vetor u)=\int_{\Omega} X(\vetor x)^T\cdot A(\vetor x)X(\vetor x)\,d\vetor x,$$
where, following the notation introduced in the Appendix,
$$H(\vetor x)=\left[\matrix{  {\partial u_1\over \partial x_1} & {\partial u_1\over \partial x_2} \cr
                              {\partial u_2\over \partial x_1} & {\partial u_2\over \partial x_2} \cr }\right]
			:=\left[\matrix {a & b+d \cr b-d & c \cr}\right]$$
and
$$X(\vetor x)^T:=(a,c,b,d)=\left({\partial u_1\over \partial x_1}, {\partial u_2\over \partial x_2},
                    {1\over 2}\left[{\partial u_1\over \partial x_2}+{\partial u_2\over \partial x_1}\right],
					{1\over 2}\left[{\partial u_1\over \partial x_2}-{\partial u_2\over \partial x_1}\right]\right)$$
			
We notice that if $A({\vetor x})-\alpha I$ is uniformly semipositive in  $\Omega$, then
$$X(\vetor x)^T\cdot A(\vetor x)X(\vetor x)\ge \alpha\left(\Mod{{\partial u_1\over \partial x_1}}^2+
    \Mod{{\partial u_2\over \partial x_2}}^2+{1\over 2}\Mod{{\partial u_1\over \partial x_2}}^2+
	   {1\over 2}\Mod{{\partial u_2\over \partial x_1}}^2\right),$$
and consequently,
$${\cal L}(\vetor u,\vetor u)\ge {\alpha\over 2}\left(\norma{\nabla u_1}{L^2}^2+\norma{\nabla u_2}{L^2}^2\right)
    ={\alpha\over 2}\norma{\vetor u}{\cal V}^2.$$

In order to analyze the matrix \cite{MatrizA} and in view of the conditions \cite{CondS1S2}, we must distinguish two cases:
$s_2<0<s_1$ and $0\le s_2<s_1$. In both cases, we fix a constant $k>\max\{0, s_2s_1^{-1}\}$ and take $\eps:=s_1-s_2k^{-1}$.
Now, let $a_0=a_0(\vetor x)$ and $b_0=b_0(\vetor x)$ be the functions defined by
$$\left\{\eqalign{
a_0 & := -2s_2\tra B_0^{-1}- 2\Bigl(s_1\sqrt{\det B_0}-s_2\sqrt{\det B_0^{-1}}\Bigr),\cr
b_0 & := -2s_2\tra B_0^{-1}+ 2\Bigl(s_1\sqrt{\det B_0}-s_2\sqrt{\det B_0^{-1}}\Bigr).\cr
}\right.\reqlabel{Defa2b2}$$

Assuming that $\det B_0\ge k$, we have the following inequalities
$$\left\{\eqalign{
b_0-a_0 &  \ge 4\sqrt{-s_1s_2}\quad\hbox{if  $s_2<0$}, \cr
b_0-a_0 & = 4\left(s_1-{s_2\over \det B_0}\right)\sqrt{\det B_0}\ge 4\eps\sqrt{k}\quad\hbox{if $s_2\ge 0$}.\cr
}\right.\reqlabel{a2-b2}$$
Notice that the above inequalities indicate that the interval $[a_0(\vetor x), b_0(\vetor x)]$ has nonempty interior
for all $\vetor x\in\Omega$ if $\det B_0\ge k$, and this will be essential in proving the next theorem.

Finally, we set
$$\overline{d}=\sup_{{\vetor x}\in\Omega}\left({\tra B_0 \over \sqrt{\det B_0}}\right).$$

\noindent{\bf Theorem\ \lemlabel{Thm1}:} {\sl Suppose that $\det B_0\ge k$. Let $\alpha>0$ such that
$$\alpha\overline{d}<\cases{2\sqrt{-\mathstrut s_1s_2} & if $s_2<0$,\cr 2\eps\sqrt{k} & if $s_2\ge 0$,\cr}\reqlabel{DefC-eps}$$
and assume that $p_0$ satisfies the condition
$$a_0({\vetor x})+\alpha\overline{d}<-2p_0({\vetor x})<b_0({\vetor x})-\alpha\overline{d},
        \quad \forall\,{\vetor x}\in\Omega.\reqlabel{CondGap1}$$
Then, there exists a constant $\beta_0=\beta_0(s_1,s_2,\alpha)\ge 0$ such that the matrix
$A({\vetor x})-\alpha I$ is uniformly positive semidefinite in $\Omega$, provided that
$\beta(\vetor x,\tau)\ge\beta_0$ for almost all $\vetor x\in\Omega$. }

\smallskip
\noindent{\bf Remark\ \lemlabel{Obs2}:} Since $2\sqrt{\det B_0}\le \tra B_0$, it follows that $\overline{d}\ge 2$. Hence,
if $\alpha$ satisfies \cite{DefC-eps}, we have necessarily
$$\alpha<\cases{\sqrt{-\mathstrut s_1s_2} & if $s_2<0$,\cr \eps\sqrt{k} & if $s_2\ge 0$.\cr}\reqlabel{DefC-eps-Bi}$$

\smallskip\noindent
{\bf Proof of Theorem~\cite{Thm1}:}  The nonzero entries of the matrix $A$ are
$$\left\{\eqalign{
A_{11} & := \beta +2s_1\gamma_1 -2s_2\gamma_1^{-1}, \cr
A_{22} & := \beta +2s_1\gamma_2 -2s_2\gamma_2^{-1}, \cr
A_{12} & := \beta +{1\over 2}\tra T_0, \cr}\right.
\qquad
\left\{\eqalign{
A_{33} & := 2s_1\tra B_0 -2s_2\tra B_0^{-1}-\tra T_0, \cr
A_{44} & := \tra T_0, \cr
A_{34} & := s_1(\gamma_2-\gamma_1)-s_2(\gamma_1^{-1}-\gamma_2^{-1}), \cr}\right.
$$			
where $\gamma_1$ e $\gamma_2$ are the eigenvalues of $B_0$. To simplify the notation, we introduce the functions
$f,g:(0,+\infty)\rightarrow\R$, as
$$f(\gamma):=s_1\gamma-s_2\gamma^{-1}\quad\hbox{\rm and}\quad g(\gamma):=s_1\gamma+s_2\gamma^{-1}.$$

A necessary and sufficient condition for the matrix $A-\alpha I$ be positive semidefinite is
$$\min\Bigl\{A_{11}-\alpha,\,\, (A_{11}-\alpha)(A_{22}-\alpha)-A_{12}^2,\,\, A_{33}-\alpha,\,\,
(A_{33}-\alpha)(A_{44}-\alpha)-A_{34}^2\Bigr\}\ge 0,\,\,\forall{\vetor x}\in\Omega.\reqlabel{Sylvester} $$
It is clear that the condition \cite{Sylvester} implies, in particular, $A_{22}-\alpha\ge 0$ and $A_{44}-\alpha\ge 0$,
since  $A$ is symmetric.

\smallskip
\noindent{\bf Step~1:} {\sl Analysis of the first block of $ A$:}

In the case $s_2<0$, we have $f(\gamma)\ge \sqrt{-\mathstrut s_1s_2}$ for all $\gamma>0$. So,
$$A_{11}-\alpha=\beta-\alpha + 2f(\gamma_1)\ge \beta-\alpha+2\sqrt{-\mathstrut s_1s_2}>\beta-\alpha.$$
In the case $s_2\ge 0$, we can assume without loss of generality that $\gamma_1\ge \gamma_2$.
Then, as $s_1-s_2/\det B_0\ge k$, we have
$$f(\gamma_1)\ge s_1\gamma_1-{s_2\over \gamma_2}\ge k\gamma_1>0,$$
and
$$A_{11}-\alpha=\beta-\alpha + 2f(\gamma_1)\ge \beta-\alpha+2k\gamma_1\ge\beta-\alpha.$$
Therefore, in the two cases, $A_{11}-\alpha \ge 0$ if $\beta\ge \alpha$.

On the other hand, if we denote $f_i=f(\gamma_i)$ and $g_i=g(\gamma_i)$, $i=1,2$, we get
$$\eqalign{
(A_{11}-\alpha)(A_{22} & -\alpha)-A_{12}^2 = (\beta-\alpha+2f_1)(\beta-\alpha+2f_2)
           -\left[\beta+{1\over 2}\tra T_0\right]^2 \cr
   & = (\beta-\alpha)\bigl[2(f_1+f_2)-2(\alpha-p_0)-(g_1+g_2)\bigr]+4f_1f_2-\left[(\alpha-p_0)+{g_1+g_2\over 2}\right]^2.\cr
}$$
Since
$$\eqalign{
2(f_1+f_2)-(g_1+g_2) & =s_1\tra B_0-3s_2\tra B_0^{-1},\cr
\alpha-p_0+{g_1+g_2\over 2} & =\alpha+{1\over 2}\tra T_0. \cr
}$$
we have
$$\eqalign{
  (A_{11}-\alpha)(A_{22}-\alpha)-A_{12}^2 & =
    (\beta-\alpha)\left[s_1\tra B_0-3s_2\tra B_0^{-1}-2(\alpha-p_0)\right]+4f_1f_2 -\left[\alpha+{1\over 2}\tra T_0\right]^2. \cr
}$$
Therefore, if
$$-2p_0<-2\alpha+s_1\tra B_0-3s_2\tra B_0^{-1},\,\,\,\forall{\vetor x}\in\Omega, \reqlabel{PrimaCond}$$
it follows that $(A_{11}-\alpha)(A_{22}-\alpha)-A_{12}^2\ge 0$ for $\beta>0$ large enough.

\smallskip\goodbreak
\noindent{\bf Step~2:} {\sl Analysis of the second block of $ A$:}

By the definition of  $T_0$, we have
$$A_{33}-\alpha = s_1\tra B_0-3s_2\tra B_0^{-1}+2p_0 -\alpha.$$
Hence, $A_{33}-\alpha \ge 0$ if, and only if,
$$-2p_0 \le -\alpha+s_1\tra B_0-3s_2\tra B_0^{-1}.\reqlabel{PrimaCondBis}$$
We notice that $A_{44}-\alpha=-2p_0s_1\tra B_0+s_2\tra B_0^{-1}-\alpha$, which implies that $A_{44}-\alpha\ge 0$
if, and only if,
$$-2p_0 \ge \alpha -s_1\tra B_0-s_2\tra B_0^{-1}.\reqlabel{SecCondBis}$$
The conditions \cite{PrimaCondBis} and \cite{SecCondBis} can be expressed by
$$-s_1\tra B_0-s_2\tra B_0^{-1}+\alpha \le 0 \le -\alpha+s_1\tra B_0-3s_2\tra B_0^{-1},$$
It is noteworthy that \cite{PrimaCond} implies \cite{PrimaCondBis}. Moreover, the interval
$$[-s_1\tra B_0-s_2\tra B_0^{-1}+\alpha\, , \, -2\alpha+s_1\tra B_0-3s_2\tra B_0^{-1}]\reqlabel{DefInter}$$
is not empty, because if we denote
$$\left\{\eqalign{
a_* & := -s_1\tra B_0-s_2\tra B_0^{-1},\cr
b_* & := s_1\tra B_0-3s_2\tra B_0^{-1},\cr
}\right.$$
it follows that
$$b_*-a_*-3\alpha=2s_1\tra B_0-2s_2\tra B_0^{-1}-3\alpha.$$
and hence, from \cite{DefC-eps-Bi}:

$$\eqalign{
\hbox{1) in the case $s_2<0$ we have\ } &  b_*-a_*-3\alpha=2f_1 + 2f_2 -3\alpha\ge 4\sqrt{-s_1s_2}-3\alpha>0,\cr
\hbox{2) in the case $s_2\ge 0$ we have\ }& b_*-a_*-3\alpha=2\tra B_0\left(s_1-{s_2\over \det B_0}\right)-3\alpha\ge 4\eps\sqrt{k}-3\alpha>0,\cr
}$$
which implies that the interval defined by \cite{DefInter} is not empty.

Now, according to the notation introduced above, we have
$$\left\{\eqalign{
A_{33} & := f_1+f_2-2s_2\tra B_0^{-1}+2p_0 \cr
A_{44} & := g_1+g_2-2p_0  \cr
A_{34} & := g_2-g_1 \cr}\right.
$$	
So,
$$\eqalign{
(A_{33}-\alpha)(A_{44}-\alpha)-A_{34}^2 & = \Bigl[f_1+f_2-2s_2\tra B_0^{-1}+2p_0-\alpha\Bigr]
          \Bigl[g_1+g_2-2p_0 -\alpha\Bigr]- [g_2-g_1]^2  \cr
    & = \bigl[F+2p_0\bigr]\bigl[G-2p_0\bigr]-[g_2-g_1]^2	\cr
}$$
where we are denoting
$$F:=f_1+f_2-2s_2\tra B_0^{-1}-\alpha\quad \hbox{\rm and}\quad   G:=g_1+g_2-\alpha.$$
Hence,
$$(A_{33}-\alpha)(A_{44}-\alpha)-A_{34}^2=FG+2p_0(G-F)-(g_2-g_1)^2-4p_0^2.$$
For $X:=-2p_0$, we have $(A_{33}-\alpha)(A_{44}-\alpha)-A_{34}^2\ge 0$ if, and only if,
$$X^2-(F-G)X +(g_2-g_1)^2-FG \le 0.\reqlabel{FormQuad}$$
The above inequality holds if the discriminant of the binomial \cite{FormQuad} is positive. In fact, we have
$$\eqalign{
(F-G)^2 & -4\bigl[(g_2-g_1)^2-FG\bigr] = \bigl[f_1+f_2-2s_2\tra B_0^{-1}+g_1+g_2-2\alpha\bigr]^2-4(g_2-g_1)^2  \cr
   & = \bigl[2s_1\tra B_0-2s_2\tra B_0^{-1}-2\alpha\bigr]^2 -4\bigl[s_1(\gamma_2-\gamma_1)+s_2(\gamma_2^{-1}-\gamma_1^{-1})\bigr]^2 \cr
   & = 4\bigl[s_1(\gamma_1+\gamma_2)-s_2(\gamma_1^{-1}+\gamma_2^{-1})-\alpha\bigr]^2 - 4\bigl[s_1(\gamma_2-\gamma_1)+s_2(\gamma_2^{-1}-\gamma_1^{-1})\bigr]^2\cr
   & = 4\bigl[2s_1\gamma_2-2s_2\gamma_1^{-1}-\alpha\bigr]\bigl[2s_1\gamma_1-2s_2\gamma_2^{-1}-\alpha\bigr]\cr
   & = 16\bigl[s_1\gamma_2-s_2\gamma_1^{-1}-\alpha/2\bigr]\bigl[s_1\gamma_1-s_2\gamma_2^{-1}-\alpha/2\bigr]\cr
   & = 16\Bigl[f\Bigl(\sqrt{\det B_0}\Bigr)^2 -
            {\alpha\over 2}\left(s_1\tra B_0-s_2\tra B_0^{-1}\right)+{\alpha^2\over 4}\Bigr].\cr
}$$

Note that
$$\eqalign{
s_1\tra B_0 - s_2\tra B_0^{-1}  & = s_1\tra B_0 -s_2\left({\tra B_0\over \det B_0}\right)
                     = \tra B_0\left(s_1-{s_2\over \det B_0}\right)  \cr
								& = \left({\tra B_0\over\sqrt{\det B_0}}\right)f\bigl(\sqrt{\det B_0}\bigr).\cr }$$
To simplify the notation, consider
$$C:=f\bigl(\sqrt{\det B_0}\bigr)\quad\hbox{\rm and}\quad D:={\tra B_0\over\sqrt{\det B_0}}.$$
Then,
$$C^2-{\alpha\over 2}CD+{\alpha^2\over 4}\ge C^2-{\alpha\over 2}CD=C^2\left(1-{\alpha D\over 2C}\right).$$

Note also that, from \cite{DefC-eps} we have:
$$\eqalign{
\hbox{\rm 1)\ \ if $s_2<0$,}\quad   & \alpha D\le \alpha\overline{d}<2\sqrt{-s_1s_2}\le 2f\bigl(\sqrt{\det B_0}\bigr)=2C;\cr
\hbox{\rm 2)\ \ if $s_2\ge 0$,}\quad& \alpha D\le \alpha\overline{d}<2\eps\sqrt{k}\le
         2\sqrt{\det B_0}\left(s_1-{s_2\over\det B_0}\right)=2f\bigl(\sqrt{\det B_0}\bigr)=2C,\cr
}$$
which  implies that, in both cases,  $0<\alpha D/2C<1$. So, by calculating the roots $a_\alpha$ and $b_\alpha$ of the binomial \cite{FormQuad}, we get,
$$\left\{\eqalign{
b_\alpha & = -2s_2\tra B_0^{-1} + 2\sqrt{\left[f\Bigl(\sqrt{\det B_0}\Bigr)\right]^2-
            {\alpha\over 2}\left(s_1\tra B_0-s_2\tra B_0^{-1}\right)+{\alpha^2\over 4}}, \cr
a_\alpha & = -2s_2\tra B_0^{-1} - 2\sqrt{\left[f\Bigl(\sqrt{\det B_0}\Bigr)\right]^2-
            {\alpha\over 2} \left(s_1\tra B_0-s_2\tra B_0^{-1}\right)+{\alpha^2\over 4}}, \cr			
}\right.\reqlabel{DuasRaizes}$$
and the condition $(A_{33}-\alpha)(A_{44}-\alpha)-A_{34}^2\ge 0$ is equivalent to $a_\alpha\le -2p_0\le b_\alpha$.

We can rewrite \cite{DuasRaizes} as
$$\left\{\eqalign{
b_\alpha & = -2s_2\tra B_0^{-1} + 2\sqrt{C^2-{\alpha\over 2}DC+{\alpha^2\over 4}}, \cr
a_\alpha & = -2s_2\tra B_0^{-1} - 2\sqrt{C^2-{\alpha\over 2}DC+{\alpha^2\over 4}}, \cr	
}\right.$$
so that
$$\sqrt{C^2-{\alpha\over 2}CD+{\alpha^2\over 4}}\ge \sqrt{C^2-{\alpha\over 2}CD}
   =C\sqrt{1-{\alpha\over 2}{D\over C}}>C-{\alpha\over 2}D.$$
Therefore, from \cite{DefC-eps} and \cite{Defa2b2}, we have
$$b_\alpha\ge b_0-\alpha\overline{d}\quad\hbox{\rm and}\quad a_\alpha\le a_0+\alpha\overline{d}.$$
Notice also that, from \cite{a2-b2} and \cite{DefC-eps}, we have $b_\alpha-a_\alpha \ge b_0-a_0 - 2\alpha\overline{d}>0$.
Thus, to conclude the proof, it suffices to show that, under the hypothesis \cite{DefC-eps}, the following inequalities hold:
$$a_*+\alpha\le a_0+\alpha\overline{d}\quad\hbox{\rm and}\quad  b_0-\alpha\overline{d} \le b_*-2\alpha.$$

Indeed, first note that
$$2\sqrt{\det B_0}  \le \tra B_0\quad \hbox{\rm and}\quad 2\sqrt{\det B_0^{-1}} \le \tra B_0^{-1}.$$
Therefore, in the case $s_2<0$, we have
$$s_1\tra B_0-s_2\tra B_0^{-1} \ge 2s_1\sqrt{\det B_0}-2s_2\sqrt{\det B_0^{-1}}=2f\bigl(\sqrt{\mathstrut\det B_0}\bigr),$$
from which we conclude that
$$s_1\tra B_0-3s_2\tra B_0^{-1} \ge -2s_2\tra B_0^{-1}+2f\bigl(\sqrt{\mathstrut\det B_0}\bigr) $$
and so, $b_*\ge b_0$.
On the other hand, it is easy to see that
$$-s_1\tra B_0+s_2\tra B_0^{-1}\le -2s_1\sqrt{\det B_0}+2s_2\sqrt{\det B_0^{-1}}=-2f\bigl(\sqrt{\mathstrut\det B_0}\bigr),$$
which implies that,
$$-s_1\tra B_0-s_2\tra B_0^{-1}\le -2s_2\tra B_0^{-1}-2f\bigl(\sqrt{\mathstrut\det B_0}\bigr), $$
and so, $a_*\le a_0$.

In the case $s_2\ge 0$, we have
$$\eqalign{
b_*-2\alpha\ge b_0-\alpha\overline{d}
     & \iff \tra B_0\left(s_1-{s_2\over\det B_0}\right)-2\alpha\ge 2f\bigl(\sqrt{\det B_0}\bigr)-\alpha\overline{d}\cr
	 & \iff{\tra B_0\over\sqrt{\det B_0}}f\bigl(\sqrt{\det B_0}\bigr)-2\alpha \ge 2f\bigl(\sqrt{\det B_0}\bigr)-\alpha\overline{d}\cr
	 & \iff \left({\tra B_0\over\sqrt{\det B_0}}-2\right)f\bigl(\sqrt{\det B_0}\bigr) \ge (2-\overline{d})\alpha.\cr		
}$$
Since $\tra B_0\ge 2\sqrt{\det B_0}$ and $\overline{d}\ge 2$, it follows that $b_*-2\alpha\ge b_0-\alpha\overline{d}$  holds for all $\alpha>0$.

Likewise,
$$\eqalign{
a_*+\alpha\le a_0+\alpha\overline{d}
     & \iff -s_1\tra B_0+s_2\tra B_0^{-1}+\alpha \le -2f\bigl(\sqrt{\det B_0}\bigr)+\alpha\overline{d}\cr
	 & \iff s_1\Bigl(2\sqrt{\det B_0}-\tra B_0\Bigr)-s_2\left({2\over\sqrt{\det B_0}}-\tra B_0^{-1}\right) \le  (\overline{d}-1)\alpha\cr
	 & \iff \left(s_1-{s_2\over\det B_0}\right)\Bigl(2\sqrt{\det B_0}-\tra B_0\Bigr) \le (\overline{d}-1)\alpha\cr		
}$$
and $a_*+\alpha\le a_0+\alpha\overline{d}$ holds for all $\alpha>0$. This finishes the proof.\quad\cqd

\smallskip
\noindent{\bf Remark\ \lemlabel{Obs3}:} The above considerations permit us to conclude (by Lax-Milgram Lemma) that the boundary value
problem \cite{ProbCont} admits a unique solution. In fact,

\smallskip
\noindent\noindent{\bf Corollary\ \lemlabel{Corol1}:} {\sl Under the hypothesis of Theorem~\cite{Thm1},
the variational problem \cite{Variacional} admits a unique solution $\vetor u\in{\cal V}$.\quad\cqd
}

\smallskip
\noindent{\bf Remark\ \lemlabel{Obs4}:} Theorem~\cite{Thm1} gives a sufficient condition for the existence of a unique weak solution of the boundary value problem
\cite{ProbCont} corresponding to each time step of the successive approximation. As we are supposing that the material is nearly incompressible,
it is reasonable to expect that $\det B_0\approx 1$. This implies that, if $\gamma_1$ and $\gamma_2$ are the eigenvalue of $B_0$,
$\gamma_1\approx 1/\gamma_2$ and $\tra B_0^{-1}\approx \gamma_1+1/\gamma_1$. So, the hypothesis \cite{CondGap1} does not means that we are assuming
that $p_0$ is small. On the other hand, numerical experiments show that the hypothesis \cite{CondGap1} can be very restrictive
in the presence of gravitational body forces. In this case, we can incorporate the potential of the gravitational force into the pressure,
and analyze the re-formulated problem.

\smallskip
\noindent{\bf Remark\ \lemlabel{Obs5}:} The previous results hold if we assume that $\Gamma_3=\emptyset$. In fact, unlike the space $\cal V$
introduced in \cite{DefV2D}, we must consider
$${\cal V}=\bigl\{\vetor u\in(H^1(\Omega))^2\,;\, \vetor u\cdot\vetor n_\kappa=0 \hbox{\ on\ } \Gamma_2\bigr\}.\reqlabel{DefVbis}$$
However, in this case, it is necessary to assume that the domain $\Omega$ satisfies a geometric property to ensure that \cite{DefNormaV} is a
norm. This can be done by supposing that $\Omega$ has the following property:
{\sl There is no constant vector $\vetor c \in \R^2$ such that $\vetor c\cdot\vetor n_\kappa(\vetor x)=0,\,\, \forall\vetor x\in\Gamma_2$.}

\bigskip\goodbreak
\noindent{\bigbf\numsection.\ Appendix}
\bigskip

Without loss of generality, we can assume that $B_0$ is a diagonal matrix, given by
$$B_0=\pmatrix{\gamma_1 & 0 \cr 0 & \gamma_2}$$
and, in this case,
$$T_0=\pmatrix{ t_1 & 0 \cr 0 & t_2}=\pmatrix{ -p_0 +f(\gamma_1)& 0 \cr 0 & -p_0+f(\gamma_2) \cr},$$
where  $f(\gamma)=s_1\gamma- s_2\gamma^{-1}$.

Writing the quadratic form \cite{DefFormQuadA} as
$${\cal A}(\vetor x, H,H) = {\cal A}_1(\vetor x,H,H)+{\cal A}_2(\vetor x,H,H)+{\cal A}_3(\vetor x,H,H)+{\cal A}_4(\vetor x,H,H),$$
where
$$\left\{\eqalign{
{\cal A}_1(\vetor x, H,H) & = \tra(H)\tra\Bigl[(T_0+\beta I)H^T\Bigr],\cr
{\cal A}_2(\vetor x, H,H) & = -\tra(T_0 H^TH^T),\cr
{\cal A}_3(\vetor x, H,H) & = s_1\tra\Bigl[(HB_0+B_0H^T)H^T\Bigr],\cr
{\cal A}_4(\vetor x, H,H) & = -{s_2}\tra\Bigl[(B_0^{-1}H+H^TB_0^{-1})H^T\Bigr].\cr
}\right.$$
and the matrix $H$ as $H=E+R$, where $E={1\over 2}(H+H^T)$ and $R={1\over 2}(H-H^T)$, with
$$E=\pmatrix{a & b \cr b & c \cr},\qquad R=\pmatrix{ 0 & d \cr -d & 0 \cr}, $$
we obtain,

\noindent 1)\quad $\tra(T_0H^T+\beta H^T)  = \tra(T_0E+\beta E) = at_1+ct_2+\beta(a+c)$, which gives
$${\cal A}_1(\vetor x, H,H)=(a+c)(at_1+ct_2)+\beta(a+c)^2.\reqlabel{E1}$$

\noindent 2)\quad Since $T_0H^TH^T=T_0(E^2+R^2)-T_0(ER+RE)$, we have $\tra(T_0H^TH^T)=\tra\bigl[T_0(E^2+R^2)\bigr]$ and
a direct calculation gives
$${\cal A}_2(\vetor x, H,H)=-t_1(a^2+b^2-d^2)-t_2(b^2+c^2-d^2).\reqlabel{E2}$$

\noindent 3)\quad We notice that
$$\left\{\eqalign{
B_0H^TH^T & = B_0(E^2+R^2)-B_0(ER+RE),\cr
HB_0H^T   & = (EB_0E-RB_0R)+(RB_0E-EB_0R).\cr
}\right.$$
Since $B_0(ER+RE)$ and $RB_0E-EB_0R$ are skew symmetric, we have
$$\tra\Bigl[(HB_0+B_0H^T)H^T\Bigr]=\tra\bigl[B_0(E^2+R^2)+(EB_0E-RB_0R)\bigr],$$
and a direct calculation gives
$${\cal A}_3(\vetor x, H,H) = 2s_1\Bigl[\gamma_1 a^2+\gamma_2 c^2+(\gamma_1+\gamma_2)b^2+(\gamma_2-\gamma_1)bd\Bigr].\reqlabel{E3}$$

\noindent 4)\quad As before,
$$\left\{\eqalign{
B_0^{-1}HH^T & = B_0^{-1}(E^2-R^2)+B_0(RE-ER+RE),\cr
H^TB_0H^T   & = (EB_0^{-1}E+RB_0^{-1}R)-(EB_0^{-1}R+RB_0^{-1}E),\cr
}\right.$$
a direct calculation gives
$${\cal A}_4(\vetor x, H,H)=-2{s_2}\bigl[\gamma_1^{-1}a^2+\gamma_2^{-1}c^2
      +(\gamma_1^{-1}+\gamma_2^{-1})b^2+(\gamma_1^{-1}-\gamma_2^{-1})bd\bigr].\reqlabel{E4}$$

Therefore, by denoting $X=(a,c,b,d)^T$ and considering \cite{E1}-\cite{E4}, we can express the quadratic form 	\cite{DefFormQuadA}
as
$${\cal A}(\vetor x, H,H)=X^T\cdot A(\vetor x)X,$$
where $A(\vetor x)$ is the matrix
$$A(\vetor x)\kern-2pt=\kern-2pt\pmatrix{
\stl \beta+2s_1\gamma_1-2{s_2}\gamma_1^{-1} & \stl \beta+{1\over 2}\tra T_0 & \stl 0 &\stl 0 \cr
\stl \beta+{1\over 2}\tra T_0 & \stl \beta+2s_1\gamma_2-2{s_2}\gamma_2^{-1}    & \stl 0 & \stl 0 \cr
\stl 0 & \stl 0 & \stl  2s_1\tra B_0 - 2{s_2}\tra B_0^{-1}-\tra T_0    & \stl  s_1(\gamma_2-\gamma_1)-{s_2}(\gamma_1^{-1}-\gamma_2^{-1})\cr
\stl 0 & \stl 0 & \stl s_1(\gamma_2-\gamma_1)-{s_2}(\gamma_1^{-1}-\gamma_2^{-1}) & \stl\tra T_0                              \cr
}$$

\bigskip\goodbreak
\noindent{\bigbf\ Acknowledgements.} The authors (ISL and MAR) are partially supported by CNPq-Brasil.
\bigskip

\bigskip\goodbreak
\noindent{\bigbf\  References}\par
\kern-.6truecm
\MakeBibliography{}

%%%  ******************************************************
%%%                    END DOCUMENT 
%%%  ******************************************************

\bye